\documentclass[preprint,11pt]{elsarticle}
\usepackage{amssymb}
\usepackage{amsmath}
\usepackage{amsfonts,graphicx,color}
\newtheorem{theorem}{Theorem}[section]
\newtheorem{lemma}{Lemma}[section]
\newtheorem{corollary}{Corollary}[section]

\newcommand{\ignore}[1]{}{}
\newproof{pf}{Proof}
\newproof{pot1}{Proof of Theorem \ref{Thm1}}
\newproof{pot2}{Proof of Theorem \ref{prop}}
\newproof{pot3}{Proof of Theorem \ref{Z-distance}}
\newtheorem{condition}{Condition}

\newcommand\beq{\begin{equation}}
	\newcommand\eeq{\end{equation}}

\makeatletter

\newcommand{\Rmnum}[1]{\expandafter\@slowromancap\romannumeral #1@}

\def\cal#1{\mathcal{#1}}
\def\geq{\geqslant}
\def\leq{\leqslant}
\def\P{\mathbb{P}}
\def\E{\mathbb{E}}

\def\R{\mathbb{R}}
\def\Z{\mathbb{Z}}

\begin{document}
	\let\today\relax
	\date{\quad}
	
	\begin{frontmatter}

	\title{From law of the iterated logarithm to Zolotarev distance for supercritical branching processes in random environment}
	\author{Yinna Ye\corref{c1}}
	\ead{yinna.ye@xjtlu.edu.cn}
	\cortext[c1]{Corresponding author}
	\address{Department of Applied Mathematics, School of Mathematics and Physics, Xi'an Jiaotong-Liverpool University, Suzhou $215123$, P. R. China}

\begin{abstract}
		Consider $(Z_n)_{n\geq0}$ a supercritical branching process in an independent and identically distributed environment. Based on some recent development in martingale limit theory, we established law of the iterated logarithm, strong law of large numbers, invariance principle and optimal convergence rate in the central limit theorem under Zolotarev and Wasserstein distances of order $p\in(0,2]$ for the process $(\log Z_n)_{n\geq0}$.
		\end{abstract}

\begin{keyword}
	Branching processes in random environment \sep Law of the iterated logarithm \sep Law of large numbers \sep Convergence rates in central limit theorem \sep Zolotarev distance \sep Wasserstein distance
\vspace{0.3cm}
\MSC 60J80; 60K37; 60F05; 62E20  
\end{keyword}

\end{frontmatter}


\section{Introduction}

Branching process in random environment (BPRE), initially introduced by \citet{SW69}, is a generalization of Galton-Watson process. The asymptotic behaviours for BPRE have been studied widely and intensively since the 1970s. For instance, one may refer to \citet{AK71b,AK71a}, \citet{T77,T88}, \citet{GL01} and \citet{BGK05} for some classical results on branching processes in independent and identically distributed (i.i.d.) environment or stationary and ergodic environment; \citet{LY10} and \citet{GLL19,GLL22} for branching processes in Markovian environment. Regarding to the recent advances especially for supercritical branching processes in i.i.d. environment, for instance, \citet{HL12,HL14} and \citet{LLGW14} established the large and moderate deviation principles, convergence theorems in $\mathbb{L}^p$ space and the central limit theorem (CLT); \citet{GLM23} studied asymptotic distribution and harmonic moments.

The main objective of this paper is to establish some additional limit theorems for supercritical BPRE, which have rarely appeared in the literature so far, including law of the iterated logarithm (LIL), strong law of large numbers (LLN), invariance principle and the rate of convergence in the CLT under Zolotarev and Wasserstein distances. The study of Wasserstein and Zolotarev distances are of independent interests, as they have been widely used in applications nowadays, such as machine learning (\cite{GAA17}), computer vision (\cite{NBCT09}, \cite{SDC15}) and microbiome studies (\cite{WCL21}). The methods used in this paper are mainly based on some recent development in martingale limit theory, especially on the invariance principles and LIL for martingales. One may refer to \citet{HH80} and \citet{S93} for the martingale limit theory and invariance principles; \citet{WMSF20} for some recent advances on a Berry-Esseen bound for martingales.

Let's firstly introduce the model of BPRE. Suppose  $\xi=(\xi_0,\,\xi_1,\,\ldots)$ is a sequence of i.i.d. random variables. Usually, $\xi$ is called environment process and  $\xi_n$ represents the random environment in the $n$-th generation. A discrete-time random process $(Z_{n})_{n\geq 0}$ is called BPRE $\xi$, if it satisfies the following recursive relation:
$$
Z_{0}=1,~Z_{n+1}=\sum_{i=1}^{Z_{n}} X_{n, i},\ \ n \geq 0,
$$
where $X_{n, i}$ represents the number of offspring produced by the $i$-th particle in the $n$-th generation.  The distribution of $X_{n,i}$ depending on the environment $\xi_{n}$ is denoted by 
$p(\xi_{n})=\{p_{k}(\xi_{n})=\mathbb{P}(X_{n,i}=k|\xi_{n}): k \in \mathbb{N}\}.$
Suppose that given $ \xi_{n} $, $ (X_{n, i})_{i\geq1}$ is a sequence of i.i.d. random variables; moreover, $ (X_{n, i})_{i\geq1}$ is independent of $(Z_{1},\ldots ,Z_{n}) $. Let $\left(\Gamma , \mathbb{P}_{\xi}\right)$ be the probability space under which the process is defined when the environment $\xi$ is given. The state space of the random environment $ \xi $ is denoted by $ \Theta $ and the total probability space can be regarded as the product space $\left(\Theta^{\mathbb{N}}\times\Gamma, \mathbb{P}\right),$ where $\mathbb{P}(\mbox{d} x, \mbox{d} \xi)=\mathbb{P}_{\xi}(\mbox{d} x) \tau(\mbox{d} \xi).$ That is, for any measurable positive function $g $ defined on $\Theta^{\mathbb{N}}\times\Gamma$, we have
$$\int g(x, \xi)\, \mathbb{P}(\mbox{d} x, \mbox{d} \xi)=\iint g(x, \xi)\, \mathbb{P}_{\xi}(\mbox{d} x)\,\tau(\mbox{d} \xi),$$
where $ \tau $ represents the distribution law of the random environment $ \xi $. And $\mathbb{P}_{\xi} $ can be regarded as the conditional probability of $\mathbb{P} $ given the environment $\xi$. The expectations with respect to (w.r.t) $ \mathbb{P}_{\xi} $ and $ \mathbb{P} $ are denoted by $ \mathbb{E}_{\xi} $ and $\mathbb{E}$, respectively. For any environment $\xi$, integer $n\geq0$ and real number $p>0$, define
$$m_n^{(p)}=m_n^{(p)}(\xi)=\sum_{i=0}^{\infty}i^p\,p_i(\xi_n),\quad m_n=m_n(\xi)=m_n^{(1)}(\xi).$$
Then
$$
m_0^{(p)}=\mathbb{E}_{\xi} Z_1^p,~~m_n=\mathbb{E}_{\xi}X_{n, i},~~i\geq 1.
$$
Consider the following random variables
$$
\Pi_{0}=1,\ \Pi_{n}=\Pi_{n}(\xi)=\prod_{i=0}^{n-1} m_{i},~~n\geq 1.
$$
It is easy to see that
$$
\mathbb{E}_{\xi}Z_{n+1}=\mathbb{E}_{\xi}\Big[\sum_{i=1}^{Z_{n}} X_{n,i} \Big]=\mathbb{E}_{\xi}\Big[\mathbb{E}_{\xi}\Big(\sum_{i=1}^{Z_{n}} X_{n,i} \Big|Z_{n}\Big)\Big]=\mathbb{E}_{\xi}\Big[\sum_{i=1}^{Z_{n}} m_{n} \Big]
=m_{n}\,\mathbb{E}_{\xi}Z_{n}.
$$
Then by recursion, we get $\Pi_{n}= \mathbb{E}_{\xi}Z_{n}$. Denote
$$
X=\log m_{0},~\mu=\mathbb{E}X,~~\sigma^2=\mathbb{E}(X-\mu)^2.
$$
The branching process $(Z_{n})_{n\geq0} $ is called supercritical, critical or subcritical according to $ \mu > 0 $, $ \mu = 0 $ or $  \mu < 0 $, respectively. Over all the paper, assume that 
\begin{equation}\label{1.1}
	p_0(\xi_0)=0\text{ a.s.}
\end{equation}
The condition (1.1) means that each particle has at least one offspring, which implies $X\geq0$ a.s. and hence $Z_n\geq1$ a.s. for any $n\geq1$.

Let's introduce now Zolotarev distance and Wasserstein distance respectively between two probability laws. For any $p>0$, denote by $[p]$ the largest integer which is strictly less than $p$, i.e.
$$[p]=\max\{n\in\Z:\, n<p\}.$$
Set $l=[p]$.  Then $p$ can be decomposed uniquely  as $p=l+\delta$, with $\delta\in(0,1]$. In the sequel, we will use this notation for the unique decomposition of any $p>0$. For a given $p>0$, consider $\Lambda_{p}$, the class of $l$-times continuously differentiable real-valued functions, defined by
$$\Lambda_p=\{f:\R\rightarrow\R:\, |f^{(l)}(x)-f^{(l)}( y)|\leq |x-y|^{p-l},\;\forall (x, y) \in \R^2\}.$$
For $l$-times continuously differentiable function $f$, set
$$\|f\|_{\Lambda_p}=\sup\left\{\frac{|f^{(l)}(x)-f^{(l)}(y)|}{|x-y|^{p-l}}:\;(x,y)\in\R^2\right\}.$$
Then for any $f\in\Lambda_p$,
\begin{align}\label{norm}
	\|f\|_{\Lambda_p}\leq1.
\end{align}
Consider a set of probability laws $\cal L(\mu, \nu)$ on $\mathbb{R}^2$ with marginals $\mu$ and  $\nu$. Zolotarev distance of order $p$ between $\mu$ and $\nu$ is defined by
$$\zeta_p(\mu,\nu)=\sup\left\{\int fd\mu-\int fd\nu: \,f\in\Lambda_p\right\}.$$
For any $p > 0$, the Wasserstein distance of order $p$ between $\mu$ and  $\nu$ is defined by 
$$W_p(\mu,\nu)=\begin{cases}\inf\left\{\int|x-y|^p\,\P(dx,dy):\,\P\in\cal L(\mu,\nu)\right\},\quad&\text{if }0<p\leq1;\\
	\inf\left\{\left(\int|x-y|^p\,\P(dx,dy)\right)^{1/p}:\,\P\in\cal L(\mu,\nu)\right\},\quad&\text{if }p>1.	
\end{cases}$$
Let $\cal F_p$ be the collection of all Borel probability measures on $\R$ with finite absolute moments of order $p$; then $(\cal F_p, W_p)$ forms a metric space, which is closely related
to the topology of weak convergence of probability distributions on $\R$.

According to \citet{KR58}, for any $0< p\leq1$,
\begin{equation}\label{eqKR}
	W_p(\mu,\nu)=\zeta_p(\mu,\nu).
\end{equation}
From \citet{R98}, it is proved that for any $p>1$,
\begin{align}\label{star}
	W_p(\mu,\nu)\leq c_p \left(\zeta_p(\mu,\nu)\right)^{1/p},
\end{align}
where $c_p>0$ is a constant depending only on $p$.

Throughout the paper, $C$ denotes a positive constant whose value may vary from line to line. If $X$ and $Y$ are two random variables, such that they follow the probability laws $P_{X}$ and $P_{Y}$ respectively, then $\zeta_p(P_{X},\,P_{Y})$ (resp. $W_p(P_X,\,P_Y)$) is simply denoted by $\zeta_p(X,\,Y)$ (resp. $W_p(X,\,Y)$). We will use $G_t$ to represent the $\mathcal{N}(0,t)$-distributed random variable with variance $t>0$; in particular, $\mathcal{N}$ represents the standard normal random variable. And $a\vee b=\max\{a,b\}$, for any $(a,b)\in\R^2$.

\section{Main results}\label{main}
In the sequel, denote by $X_i=\log m_{i},\ i\geq 0.$ Evidently, $(X_i)_{i\geq 0}$ is a sequence of i.i.d. random variables depending only on the environment $ \xi $. Let $(S_n)_{n\geq 0}$ be the random walk associated with the branching process, which is defined as follows:
$$
S_{0}=0,\quad S_{n}=\log \Pi_{n}=\sum_{i=0}^{n-1}X_i,~~n\geq 1.
$$
Then we have the following decomposition of $\log Z_n$:
\begin{equation}\label{2.1}
	\log Z_{n}=S_n+\log W_{n},
\end{equation}
where $W_{n}=\frac{Z_{n}}{\Pi_{n}}. $ The normalized population size $ (W_{n})_{n\geq0} $ is a non-negative martingale under both $  \mathbb{P} $ and $ \mathbb{P}_{\xi} $, w.r.t the natural filtration $(\cal{F}_n)_{n\geq0}$, defined by
$$
\cal{F}_{0}=\sigma\{\xi\}, \ \cal{F}_{n}=\sigma\{\xi, X_{k, i}, 0 \leq k \leq n-1, i \geq 1\},~~n \geq 1.
$$
By Doob's martingale convergence theorem and Fatou's lemma, we can obtain that $ W_{n} $ converges a.s. to a finite limit $ W $ and $\mathbb{E}W \leq 1 $. We assume the following conditions throughout this paper:
\begin{equation}\label{2.2}
	\sigma\in(0,\infty)~~\text{and}~~\mathbb{E} \frac{Z_{1}}{m_{0}} \log Z_{1}<\infty.
\end{equation}
The first condition above together with (\ref{1.1}) imply in particular that
$$
Z_1\geq1~~\text{a.s.}~~\text{and}~~\mathbb{P}(Z_1=1)=\mathbb{E}p_1(\xi_0)<1.
$$
The second condition in (\ref{2.2}) implies that  $W_{n}$ converges to $W$ in $\mathbb{L}^{1}$ and
$$
\mathbb{P}(W>0)=\mathbb{P}(Z_{n}\stackrel{n\rightarrow\infty}{\longrightarrow}\infty)=\lim_{n\rightarrow\infty}\mathbb{P}(Z_n>0)>0
$$
(See e.g. \citet{T88}). Therefore, it follows with the assumption (\ref{1.1}) that $W>0$ and $Z_{n}\stackrel{n\rightarrow\infty}{\longrightarrow}\infty$ a.s.

In the sequel, we will need the following conditions.
\begin{condition}\label{Con1}
	There exists a constant  $\delta \in (0,1)$, such that
	$$
	\mathbb{E}X^{2+\delta}=\mathbb{E}\left(\log m_0\right)^{2+\delta}<\infty.
	$$
\end{condition}
\begin{condition}\label{Con2}
	There exist some constants $p>1$ and $c>0$ such that
	$$
	\E m_0^c<\infty\quad\text{and}\quad\mathbb{E} \left(\frac{m_0^{(p)}}{m_{0}^p} \right)^c<\infty.$$
\end{condition}
We have the following LIL for the process $(\log Z_n)_{n\geq0}$.
\begin{theorem}\label{Thm1}
	Suppose that Condition \ref{Con1} is satisfied. Then 
	$$\limsup_{n\rightarrow\infty}\frac{\log Z_n-n\mu}{\sqrt{n\log\log n}}=+\sqrt{2}\,\sigma\quad\text{a.s.}$$
	and
	$$\liminf_{n\rightarrow\infty}\frac{\log Z_n-n\mu}{\sqrt{n\log\log n}}=-\sqrt{2}\,\sigma\quad\text{a.s.}$$
\end{theorem}
From the LIL theorem above, we can find immediately the following strong LLN for $(\log Z_n)_{n\geq0}$.
\begin{corollary}
	Under the same condition of Theorem \ref{Thm1}, then
	$$\lim_{n\rightarrow\infty}\frac{\log Z_n}{n}=\mu\quad\text{a.s.}$$
\end{corollary}

We also have the following invariance principle for the process $(\log Z_n-n\mu)_{n\geq0}$, which only requires the existence of the second order moment $\E X^2<\infty$. And this has been insured in the assumption (\ref{2.2}).  
\begin{theorem}\label{prop}
	Let $(B(t))_{t\geq0}$ be a standard Brownian motion. Then there exists a common probability space for $(\log Z_n)_{n\geq0}$ and $(B(t))_{t\geq0}$, such that the following holds: 
	$$\log Z_n-n\mu-B(n\sigma^2)=o\left(\sqrt{n\log\log n}\right)\quad\text{a.s.,}$$
	as $n\rightarrow\infty$.
\end{theorem}
We have the following convergence rates of $(\log Z_n)_{n\geq0}$ in the CLT under Zolotarev and Wasserstein distances respectively.
\begin{theorem}\label{Z-distance}
	Under Conditions \ref{Con1} and \ref{Con2}, we have for any $r\in[\delta,\,2]$, 
	$$\zeta_r\left(\frac{\log Z_n-n\mu}{\sqrt{n}\sigma},\,\mathcal{N}\right)\leq \frac{C}{n^{\delta/2}}.$$
\end{theorem}

\begin{corollary}
	Under Conditions \ref{Con1} and \ref{Con2}, then
	$$W_r\left(\frac{\log Z_n-n\mu}{\sqrt{n}\sigma},\;\mathcal{N}\right)\leq\frac{C}{n^{\delta/2}},\quad\text{for }r\in [\delta,\,1];$$
	$$W_r\left(\frac{\log Z_n-n\mu}{\sqrt{n}\sigma},\;\mathcal{N}\right)\leq\frac{C}{n^{\delta/2r}},\quad\text{for }r\in (1,\,2].$$
\end{corollary}
The result above is an immediate consequence of Theorem \ref{Z-distance}, the equality  (\ref{eqKR}) and the inequality (\ref{star}) respectively for the cases $r\in[\delta,\,1]$ and $r\in(1,\,2]$.

\section {Proof of main results}\label{proof}
Consider a sequence of martingale differences $(\varepsilon_i)_{i\geq 1}$, defined on a probability space $(\Omega,\cal F, \P)$, w.r.t a filtration $\cal G=(\cal G_i)_{i\geq0}$. Set 
$$M_0=0, \quad M_n=\sum_{i=1}^n\varepsilon_i, \text{ for } n\geq1.$$ Then $(M_n,\cal G_n)_{n\geq0}$ is a martingale. Set $V_n^2=\sum_{i=1}^n\varepsilon_i^2$ and $s_n^2=\E M_n^2$. We have the following lemma (see also Corollary 4.2 and Theorem 4.8 in \citet{HH80}), from which Theorem {\ref{Thm1}} can be obtained.
\begin{lemma}[\citet{HH80}]\label{ThmLIL} Suppose that $(M_n,\,\cal G_n)_{n\geq 0}$ is a zero-mean, square integrable martingale. If 
	\begin{equation}\label{eq1}
		s_n^{-2}V^2_n\xrightarrow{n\rightarrow\infty} \eta^2>0\quad\text{a.s.},
	\end{equation}
	\begin{equation}\label{eq2}
		\text{for any }\varepsilon>0,\quad \sum_{n=1}^{\infty}s_n^{-1}\E\left(|\varepsilon_n|\mathbf{1}_{\{|\varepsilon_n|>\varepsilon s_n\}}\right)<\infty
	\end{equation}
	and for some $\tau>0$,
	\begin{equation}\label{eq3}
		\sum_{n=1}^{\infty}s_n^{-4}\,\E\left(\varepsilon_n^4\,\mathbf{1}_{\{|\varepsilon_n|\leq\tau s_n\}}\right)<\infty;
	\end{equation}
	then
	$$\limsup_{n\rightarrow\infty}\left(\phi(V_n^2)\right)^{-1}M_n=+1\quad\text{a.s.}$$
	and
	$$\liminf_{n\rightarrow\infty}\left(\phi(V_n^2)\right)^{-1}M_n=-1\quad\text{a.s.},$$
	where $\phi(t)=(2t\log\log t)^{1/2}$ for any $t>0$. 	
\end{lemma}

Recall now that the random walk associated with the branching process
$$S_0=0,\quad S_{n}=\sum_{i=0}^{n-1}X_i,~~n\geq 1$$
is the sum of i.i.d. random variables $(X_i)_{i\geq0}$. In the sequel, we will use the following notations. Let $\varepsilon_i=X_{i-1}-\mu$, $i\geq 1$ and
$$M_0=0,\quad M_n=\sum_{i=1}^n\varepsilon_i,\quad n\geq1.$$
Then from (\ref{2.1}), we have

\begin{align}\label{meandecomposition}
	\log Z_n-n\mu=& M_n+\log W_n.
\end{align}
Since $(\varepsilon_i)_{i\geq1}$ is a sequence of i.i.d. zero-mean random variables and $\sigma<\infty$, $(M_n)_{n\geq0}$ becomes a zero-mean and square integrable martingale w.r.t the filtration $(\cal {F}_n)_{n\geq0}$. 

\begin{pot1} 
	Since $V_n^2=\sum_{i=1}^n(X_i-\mu)^2=\sum_{i=2}^n\varepsilon^2_i$ and $s_n^2=n\sigma^2$, from the strong LLN for $V^2_n$, 
	\begin{align}\label{LLN_V}
		s^{-2}_nV^2_n&=\frac{1}{n\sigma^2}\sum_{i=2}^n \varepsilon^2_i\xrightarrow{n\rightarrow\infty} 1\quad\text{a.s.}
	\end{align}
	the condition (\ref{eq1}) is hence satisfied, with $\eta^2=1>0$. Let $\sigma^{(2+\delta)}=\E|X-\mu|^{2+\delta}$. Then under Condition \ref{Con1}, $\sigma^{(2+\delta)}<\infty$. And the condition (\ref{eq2}) is satisfied, because for any $\varepsilon>0$,
	
	\begin{align*}
		\sum_{n=1}^{\infty}s_n^{-1}\,\E\left(|\varepsilon_n|\mathbf{1}_{\{|\varepsilon_n|>\varepsilon s_n\}}\right)\leq\frac{\sigma^{(2+\delta)}}{\varepsilon^{1+\delta}} \sum_{n=1}^{\infty} \frac{1}{s_n^{2+\delta}} =\frac{\sigma^{(2+\delta)}}{\varepsilon^{1+\delta}\,\sigma^{2+\delta}}\sum_{n=1}^{\infty} \frac{1}{n^{1+\delta/2}}<\infty.
	\end{align*}
	The condition (\ref{eq3}) is also satisfied, because for any $\tau>0$, we have
	\begin{align*}
		\sum_{n=1}^{\infty}s_n^{-4}\,\E\left(\varepsilon_n^4\,\mathbf{1}_{\{|\varepsilon_n|\leq\tau s_n\}}\right)\leq\tau^{2-\delta}\sigma^{(2+\delta)} \sum_{n=1}^{\infty}\frac{1}{s_n^{2+\delta}}=\frac{\tau^{2-\delta}\,\sigma^{(2+\delta)}}{\sigma^{2+\delta}}\sum_{n=1}^{\infty}\frac{1}{n^{1+\delta/2}}<\infty.
	\end{align*}
	
	Now applying Lemma \ref{ThmLIL} to the martingale $(M_n,\,\cal F_n)_{n\geq0}$, using the facts that $\displaystyle\lim_{n\rightarrow\infty}\frac{\phi(V_n^2)}{\phi(n)}=\lim_{n\rightarrow\infty}\sqrt{\frac{V_n^2}{n}}=\sigma$ a.s., $\displaystyle\lim_{n\rightarrow\infty}\frac{|\log W_n|}{\sqrt{n\log\log n}}=0$ a.s., we come to the statement of Theorem \ref{Thm1}. \hfill$\Box$
\end{pot1}

Theorem \ref{prop} can be proved by using Theorem 2.1 in \citet{S93}, where the following invariance principle for the martingale $(M_n,\,\cal G_n)_{n\geq 0}$ with the martingale difference sequence $(\varepsilon_i)_{i\geq1}$ is established. 

\begin{lemma}[\citet{S93}, Theorem 2.1]\label{lem3.1}
	Suppose that $\{M_n=\sum_{i=1}^n\varepsilon_i,\,\cal G_n\}_{n\geq 0}$ is a square integrable martingale and $(B(t))_{t\geq0}$ is a standard Brownian motion. Let $Z(t)=\sum_{k\leq t}\varepsilon_k$, $b(t)=\E Z^2(t)$, $t\geq0$. Assume that there exists a non-decreasing positive valued sequence $(c_n)_{n\geq1}$ with $c_n\rightarrow \infty$, such that, as $n\rightarrow\infty$,
	$$\sum_{i=1}^{n}\left[\E\left(\varepsilon_i^2|\cal G_{i-1}\right)-\E \varepsilon_i^2\right]=o(c_n)\quad\text{a.s.}$$
	and
	$$\sum_{i=1}^{\infty}c_i^{-v}\,\E|\varepsilon_i|^{2v}<\infty\text{ for some }1\leq v\leq 2.$$ 
	Then there exists a common probability space for $(\varepsilon_i)_{i\geq1}$ and $(B(t))_{t\geq0}$ such that, as $n\rightarrow\infty$,
	$$Z(n)-B(b_n)=o\left(\left(c_n\left(\log \frac{b_n}{c_n}+\log\log c_n\right)\right)^{1/2}\right)\quad\text{a.s.}
	$$
\end{lemma}

\begin{pot2}
	From (\ref{meandecomposition}), we have
	$$\log Z_n-n\mu-B(n\sigma^2)=M_n+\log W_n-B(n\sigma^2).$$
	Applying Lemma \ref{lem3.1} to the martingale $(M_n,\,\cal F_n)_{n\geq0}$, letting $c_n=n$, $b_n=n\sigma^2$ and  $v=1+\delta/2\in(1,\,3/2)$, and taking into account that $ W_{n} $ converges a.s. to a finite limit $W$, we can obtain immediately the result.
	\hfill$\Box$
\end{pot2}

To prove Theorem \ref{Z-distance}, the following lemmas will be used.
\begin{lemma}\label{lem_Yinna}
	Assume Condition \ref{Con2} hold. For the constants $p$ and $c$ given from Condition \ref{Con2}, there exists a constant $a_{p,c}\in(0,1)$ depending on $p$ and $c$, such that for any $q\in(0,\,3+a_{p,c})$,
	$$\E|\log W|^q<\infty\quad\text{and}\quad \sup_{n\geq0}\,\E|\log W_n|^q<\infty.$$
\end{lemma}
\begin{pf}
	By Lyapunov's inequality, it is enough to prove the assertion for $q\in(3,\,3+a_{p,c})$. Let $q\in(3,\,3+a_{p,c})$. Using truncation, we have
	\begin{equation*}
		\mathbb{E}|\log W|^{q}=\mathbb{E}|\log W|^{q} \mathbf{1}_{\{W>1\}}+\mathbb{E}|\log W|^{q} \mathbf{1}_{\{W\leq1\}}.
	\end{equation*}
	For the first term on the right-hand side above, we have
	\begin{equation*}
		\mathbb{E}|\log W|^{q} \mathbf{1}_{\{W>1\}}\leq   C\,  \mathbb{E}W<\infty.
	\end{equation*}
	For the second term, let $\psi$ be the annealed Laplace transform of $W$, defined by 
	$$\psi(t)=\E\mbox{e}^{-t W},\quad t\geq0.$$ 
	By Markov's inequality, we have for $t>0$,
	$$\mathbb{P}(W\leq t^{-1})\leq e\,\mathbb{E}e^{-tW}=e\,\psi(t).$$ 
	Then
	\begin{align}\label{liu1}
		\mathbb{E}|\log W|^{q}\mathbf{1}_{\{W\leq1\}}&=\, q\int_1^{\infty}\frac{1}{t}(\log t)^{q-1}\,\mathbb{P}(W\leq t^{-1})\, dt \notag \\
		&\leq\, q\, e\int_1^{\infty}\frac{\psi(t)}{t}(\log t)^{q-1}\, dt. 
	\end{align}
	From (3.20) of \citet{GLM23}, it is proved that under Condition \ref{Con2},
	\begin{align}\label{liu2}
		\psi(t)\leq Ct^{-a_{p,c}},
	\end{align}
	for some constants $C>0$, $a_{p,c}\in(0,1)$ depending on $p$ and $c$, and for all $t>0$. Therefore, combining (\ref{liu1}) with (\ref{liu2}), we get
	\begin{align*}
		\mathbb{E}|\log W|^{q}\mathbf{1}_{\{W\leq1\}}\leq\, C\int_1^{\infty}\frac{(\log t)^{q-1}}{t^{1+a_{p,c}}}\,dt
		=\,&C\int_0^{\infty}y^{q-1}\mbox{e}^{-y\,a_{p,c}}\,dy=\,Ca_{p,c}^{-q}\,\Gamma(q)<\infty,
	\end{align*}
	where $\Gamma$ is the Gamma function defined by $\Gamma(z)=\int_0^{\infty}x^{z-1}\mbox{e}^{-x}\,dx$, for $\Re{(z)}>0$. We hence obtain
	\begin{align}\label{eqliu4}
		\mathbb{E}|\log W|^{q}\mathbf{1}_{\{W\leq1\}}<\infty.
	\end{align}
	
	Next, let's prove 
	\begin{equation}\label{truncation}
		\sup_{n\geq0}\,\E|\log W_n|^{q}<\infty.
	\end{equation} 
	Notice that for any $p\geq1$, $x\mapsto|\log x|^p\,\mathbf{1}_{\{0<x\leq1\}}$ is non-negative and convex function. According to Lemma 2.1 in \cite{HL12}, we have
	$$\sup_{n\geq0}\mathbb{E}\left|\log W_{n}\right|^{q}\mathbf{1}_{\{W_{n}\leq1\}}=\mathbb{E}\left|\log W\right|^{q}\mathbf{1}_{\{W\leq1\}}.$$
	With the similar truncation as for $\mathbb{E}\left|\log W\right|^{q}$ above and by (\ref{eqliu4}), we can obtain (\ref{truncation}).
	\hfill$\Box$
\end{pf}

\begin{lemma}\label{Yinna}
	Suppose that $r\in(1,\,2]$. If $f\in\Lambda_r$, then for any $(x,y)\in\R^2$, there exists $t_0$ between $x$ and $x+y$, such that
	\begin{align*}
		|f(x+y)-f(x)|\leq C|t_0|^{r-1}|y|.
	\end{align*}
\end{lemma}
\begin{pf}
	Let $r\in(1,\,2]$, then $[r]=1$. From the definition of $\Lambda_r$ and (\ref{norm}), if $f\in\Lambda_r$, then for any $t\in\R$, $\displaystyle|f'(t)-f'(0)|\leq|t|^{r-1};$ which implies 
	\begin{align}\label{basic1}
		|f'(t)|\leq C|t|^{r-1}.
	\end{align}
	From the mean value theorem and the inequality (\ref{basic1}), for any pair of real numbers $(x,y)\in \R^2$, there exists $t_0$ between $x$ and $x+y$, such that
	$$|f(x+y)-f(x)|=|f'(t_0)||y|\leq C|t_0|^{r-1}|y|.$$
	\hfill$\Box$
\end{pf}

\begin{pot3}
	For any $n\geq0$, let $a_n=\E\log W_n$ and $R_n=\log W_n- a_n$. Then $(R_n)_{n\geq0}$ becomes a sequence of zero-mean random variables. Moreover, by (\ref{meandecomposition}), we have the following decomposition, for $n\geq0$,
	\begin{align}\label{eq4}
		\log Z_n-n\mu-a_n= M_n+R_n.
	\end{align}
	By the triangle inequality and (\ref{eq4}), 
	\begin{align}\label{eq5}
		\zeta_r(\log Z_n-n\mu,\,G_{n\sigma^2})\leq&\zeta_{r}(M_n+R_n+a_n,\,M_n+R_n)+\zeta_{r}(M_n+R_n,\,M_n)+\zeta_r(M_n,\,G_{n\sigma^2})\notag\\
		=&I_{1,r,n}+I_{2,r,n}+I_{3,r,n},
	\end{align} 
	where $I_{1,r,n}=\zeta_{r}(M_n+R_n+a_n,\,M_n+R_n)$, $I_{2,r,n}=\zeta_{r}(M_n+R_n,\,M_n)$ and $I_{3,r,n}=\zeta_r(M_n,\,G_{n\sigma^2})$. We will estimate respectively $I_{1,r,n}$, $I_{2,r,n}$ and $I_{3,r,n}$.
	
	For $I_{3,r,n}$, recall that under Condition \ref{Con1}, $(M_n)_{n\geq1}$ is a centred martingale with i.i.d. martingale differences in $\mathbb{L}^{2+\delta}$ and $\mbox{Var}(M_n)=n\sigma^2$.  According to Theorem 2.1 of \citet{DMR09}, with $p=2+\delta$, we can obtain for any $r\in[\delta,\,2+\delta]$,
	\begin{align}\label{eq6}
		\zeta_r(n^{-1/2}M_n,\,G_{\sigma^2})\leq \frac{C}{n^{\delta/2}}.
	\end{align}
	Notice that for any pair of real-valued random variables $(X,Y)$ and real constant $a$, \begin{equation}\label{scale}
		\zeta_r(aX,aY)=|a|^r\zeta_r(X,Y).
	\end{equation}
	Therefore, the inequality (\ref{eq6}) implies that for any $r\in[\delta,\,2+\delta]$, 
	\begin{align}\label{Thm3.2_0}
		I_{3,r,n}=n^{r/2}	\zeta_r(n^{-1/2}M_n,\,G_{\sigma^2})\leq Cn^{\frac{r-\delta}{2}}.
	\end{align}
	
	Let's estimate now $I_{1,r,n}$. Suppose that $r\in[\delta,\,2]$ and $f\in\Lambda_r$. We will discuss in the following two cases: when $r\in[\delta,\,1]$ or $r\in(1,\,2]$, respectively. If $r\in[\delta,\,1]$, then $[r]=0$. By the definition of $\Lambda_r$, we have for any $(x,y)\in\R^2$,
	\begin{align}\label{r-zero}
		|f(x)-f(y)|\leq|x-y|^r.
	\end{align}
	Since $M_n+R_n$ and $M_n+R_n+a_n$ are both $\mathcal{F}_n$-measurable, by using (\ref{r-zero}), we have for any $r\in[\delta,\,1]$,
	\begin{align*}
		I_{1,r,n}\leq&\sup_{f\in\Lambda_r}\E\,[\E(|f(M_n+R_n+a_n)-f(M_n+R_n)|\;|\mathcal{F}_n)]\leq|a_n|^r.
	\end{align*}
	Since for $r\in[\delta,\,1]$, the function $x\mapsto x^r$ is increasing on $(0,\infty)$, and by using Lemma \ref{lem_Yinna}, we get
	\begin{align}\label{Case1}
		I_{1,r,n}\leq|\E\log W_n|^r\leq&(\E|\log W_n|)^r\notag\\
		<&\infty.
	\end{align}
	If $r\in(1,\,2]$, then $[r]=1$. By Lemma \ref{Yinna}, we have there exists $Y_{1,n}$ between $M_n+R_n$ and $M_n+R_n+a_n$, such that
	$$\E(|f(M_n+R_n+a_n)-f(M_n+R_n)|\,|\mathcal{F}_n)\leq\, C\,|a_n|\,|Y_{1,n}|^{r-1}.$$
	Hence, 
	\begin{align}\label{Ine-I0}
		I_{1,r,n}\leq&\sup_{f\in\Lambda_r}\E|f(M_n+R_n+a_n)-f(M_n+R_n)|\notag\\
		\leq& C\,|a_n|\,\E|Y_{1,n}|^{r-1}.
	\end{align}
	Now, let's study $\E|Y_{1,n}|^{r-1}$, with $r-1\in(0,\,1]$. Since for $\delta\in(0,\,1]$, the function $x\mapsto x^{\delta}$ is increasing on $(0,\infty)$, we have
	\begin{align}
		\E|Y_{1,n}|^{r-1}\leq&\max\{\,\E|M_n+R_n+a_n|^{r-1},\,\E|M_n+R_n|^{r-1}\,\} \notag\\
		\leq& \E|M_n+R_n+a_n|^{r-1}+\E|M_n+R_n|^{r-1}.
	\end{align}
	For $r\in(1,\,2]$, applying Lyapunov's inequality to the two terms on the right-hand side, we obtain
	\begin{align}\label{Ine-I2}
		\E|M_n+R_n+a_n|^{r-1}+\E|M_n+R_n|^{r-1}\leq(\E|M_n+R_n+a_n|)^{r-1}+(\E|M_n+R_n|)^{r-1}.
	\end{align}
	For the first term above, from the assumption (\ref{2.2}), the LLN for $S_n$ and Lemma \ref{lem_Yinna}, we have
	\begin{align}\label{Ine-I3}
		(\E|M_n+R_n+a_n|)^{r-1}=(\E|S_n-n\mu+\log W_n|)^{r-1}\leq&(\E|S_n-n\mu|+\E|\log W_n|)^{r-1}\notag\\
		<&\infty.
	\end{align}
	With the same arguments, we can obtain  
	\begin{align}\label{Ine-I4}
		(\E|M_n+R_n|)^{r-1}<\infty.
	\end{align}
	Using Lemma \ref{lem_Yinna} again, we also have
	\begin{align}\label{Ine-I5}
		|a_n|=|\E\log W_n|\leq\E|\log W_n|<\infty.
	\end{align}
	Combining (\ref{Case1}) for the case $r\in[\delta,\,1]$ with (\ref{Ine-I0}) - (\ref{Ine-I5}) for the case $r\in(1,\,2]$, we can conclude that for any $r\in[\delta,\,2]$,
	\begin{align}\label{I1}
		I_{1,r,n}\leq C.
	\end{align}
	
	Now, let's deal with $I_{2,r,n}$. By Lemma \ref{lem_Yinna}, we have $\sup_{n\geq0}\|R_n\|_{2+\delta}<\infty$. Taking into account Remark 5.1 and applying Items 1 and 2 of Lemma 5.2 in \citet{DMR09} with $p=\delta+2$ therein, then we can obtain for any $r\in[\delta,\,2]$,
	\begin{align}\label{eq11}
		I_{2,r,n}=\zeta_r(M_n+R_n,\,M_n)\leq C n^{\frac{r-\delta}{2}}.
	\end{align}
	Now, combining (\ref{eq5}), (\ref{Thm3.2_0}), (\ref{I1}) and (\ref{eq11}), we have for any $r\in[\delta,\,2]$,
	\begin{align*}
		\zeta_r(\log Z_n-n\mu,\; G_{n\sigma^2})\leq Cn^{\frac{r-\delta}{2}}.
	\end{align*}
	Therefore, the theorem is derived by taking into account
	$$\zeta_r\left(\frac{\log Z_n-n\mu}{\sqrt{n}\sigma},\,\mathcal{N}\right)=\left(\frac{1}{\sqrt{n}\sigma}\right)^r\zeta_r(\log Z_n-n\mu,\, G_{n\sigma^2}).$$
	\hfill$\Box$
\end{pot3}








\end{document}